\documentclass[default,iicol]{sn-jnl}

\jyear{2022}%

\theoremstyle{thmstyleone}%
\theoremstyle{thmstyletwo}%

\theoremstyle{thmstylethree}%
\begin{document}
\title[ ]{Integral equation method for the
1D steady-state Poisson-Nernst-Planck equations} 


\author*[1]{\fnm{Zhen} \sur{Chao}}\email{zhench@umich.edu}

\author[2]{\fnm{Weihua} \sur{Geng}}\email{wgeng@smu.edu}

\author[1]{\fnm{Robert} \sur{Krasny}}\email{krasny@umich.edu}

\affil[1]{\orgdiv{Department of Mathematics}, \orgname{University of Michigan},   \orgaddress{\street{} \city{Ann Arbor}, \state{MI} \postcode{48109}, \country{USA}}}

\affil[2]{\orgdiv{Department of Mathematics}, \orgname{Southern Methodist University}, \orgaddress{\street{} \city{Dallas}, \state{TX} \postcode{75275}, \country{USA}}}


\abstract{
An integral equation method is presented for the
1D steady-state Poisson-Nernst-Planck equations
modeling ion transport through membrane channels.   
The differential equations are recast 
as integral equations using Green's 3rd~identity
yielding a fixed-point problem 
for the electric potential gradient and ion concentrations.
The integrals are discretized 
by a combination of midpoint and trapezoid rules
and
the resulting algebraic equations are solved by Gummel iteration.
Numerical tests for
electroneutral and non-electroneutral systems
demonstrate the method's 2nd order accuracy
and 
ability to resolve sharp boundary layers.
The method is applied to a
1D model of the K$^+$~ion channel with a 
fixed charge density that ensures cation selectivity.
In these tests the proposed integral equation method yields
potential and concentration profiles
in good agreement with published results.
}




\maketitle

\section{Introduction}\label{sec1}

Ion transport through membrane channels
is a fundamental process in cell
biology~\cite{hille1978ionic,hille2001ion}. 
However,
atomistic simulations of this process are expensive
and
an alternative continuum approach is based on the
Poisson-Nernst-Planck (PNP) equations,
where
the Poisson equation for the electric potential
is coupled
to the Nernst-Planck equation 
for the ion concentrations.
The continuum PNP model accounts for two effects,
(1) molecular diffusion due to ion concentration gradients,
and
(2) ion drift in the applied
and
self-consistent electric fields.
The PNP equations have been employed to study
many aspects of ion channel selectivity~\cite{liu2017incorporating},
including
the effect of fixed charge~\cite{
eisenberg_poissonnernstplanck_2007},
ion size~\cite{horng2012pnp}, 
and
membrane surface charge~\cite{xie2020finite}.  

Despite the relative efficiency in comparison with
atomistic simulations,
numerical solution of the PNP equations is still challenging
due to several factors,
(1) nonlinear coupling of the electric potential
and
ion concentrations, 
(2) the need to resolve the channel geometry,
as well as boundary layers
and
internal layers at material interfaces,
and 
(3) maintaining stability 
in convection-dominated cases~\cite{wang2021stabilized}.
Many PNP solvers have been developed using
finite-difference methods~\cite{
zheng_second-order_2011,
flavell2014conservative}, 
finite-element methods~\cite{
lu2007electrodiffusion,
lu_poissonnernstplanck_2010,
xie2020finite, 
xie2020effective}, 
and 
finite-volume methods~\cite{chainais2003convergence,
song2018electroneutral,
song2019selectivity}.
%
%
It is also important for numerical PNP solvers 
to preserve positivity of the ion concentrations
and
conserve mass in time-dependent simulations,
and
several implicit finite-difference schemes 
have been developed for this
purpose~\cite{flavell2014conservative,hu2020fully}.


The boundary element method (BEM) 
is another approach to solving linear elliptic 
boundary value problems,
where the differential equation is recast as 
an integral equation on the domain boundary
by convolution with an 
appropriate Green's function~\cite{steinbach2007numerical}.
While the BEM requires 
discretization of singular integrals,
it has the advantage of 
reduced spatial dimension.
Closely related to the PNP equations is the
Poisson-Boltzmann (PB) equation for the
electrostatic potential of a solvated protein.
The BEM has been applied to solve the 
3D linear PB equation 
on the molecular surface of a protein
(see~\cite{juffer1991electric,geng2013treecode}
and references therein),
but aside from some 
early work using a Green's function
to solve the PB equation 
for an ion channel in 1D~\cite{barcilon1992ion},
as far as we know
this approach has not been further developed
for 2D and 3D ion channel problems.

In fact there are two challenges to applying 
an integral equation approach
to the PNP system,
(1) nonlinear coupling of the
electric potential and ion concentrations,
and
(2) the ions are distributed throughout
the channel and solvent domain requiring
costly volume integrals in 2D and 3D.
Nonetheless we believe these obstacles can be overcome
using iteration to handle the nonlinearity,
and
a fast summation method with adaptive quadrature
for the volume integrals.
As a step in that direction,
the present work develops a novel integral equation method
for the 1D steady-state PNP equations.

The paper is organized as follows. 
Section~2 introduces the 1D steady-state PNP equations.
Section~3 uses Green's 3rd~identity to recast the 
PNP equations as a fixed-point problem
involving integral operators
for the potential gradient and ion concentrations. 
Section~4 presents the numerical implementation. 
Section~5 presents tests for 
electroneutral and
non-electroneutral systems with
Robin boundary conditions on the potential
and
no-flux boundary conditions on the ion concentrations.
Section~6 presents a K$^+$~ion channel simulation
following~\cite{gardner2004electrodiffusion}
and
the work is summarized in Section~7.


\section{PNP equations}\label{sec2}

This work focuses on the 
1D steady-state PNP equations 
with two ion species.
The electrostatic potential $\phi(x)$ satisfies the 
Poisson equation,
\begin{equation}\label{initialP}
-\left(\epsilon\phi^\prime(x)\right)^\prime =
z_1 e c_1(x) + z_2 e c_2(x),
\end{equation}
where $\epsilon$ is the permittivity,
the valence is
$z_1 = -1$ for anions
and
$z_2 = 1$ for cations,
$e$ is the proton charge,
$c_1(x), c_2(x)$ are the ion concentrations,
and
the domain is $[-L,L]$.
This is supplemented with the Nernst–Planck equation,
\begin{equation}\label{initialNP}
\left[D_{i}\left(c_{i}^\prime(x)+\frac{z_{i} e}{k_{B} T} c_{i}(x)\phi^\prime(x)\right)\right]^\prime=0,
\end{equation}
for $i = 1,2$, where 
$D_i, k_B, T$ are the 
diffusion coefficient of ion species $i$, 
Boltzmann constant,
and 
absolute temperature.




Following~\cite{flavell2014conservative}
the PNP equations are non-dimensionalized with respect to the
bulk concentration $c_0$ [ions/\AA$^3$], 
channel half-length $L$ [\AA],
characteristic diffusion coefficient $D_0$ [\AA$^2$/s],
characteristic potential $\phi_0$ [V],
and
water permittivity $\epsilon_t$ [F/\AA].
The resulting non-dimensional equations are
\begin{subequations}
\label{eqn:PNP_system}
\begin{align}
\label{eqn:P_non_dimensional}
& -\epsilon\phi^{\prime\prime}(x) = 
\chi_2 \left(z_1c_1(x) + z_2c_2(x) \right), \\[5pt]
\label{eqn:NP_non_dimensional} 
& \left( D_i\left( c_i^\prime(x) + 
\chi_1 z_i c_i(x) \phi^\prime(x)\right)\right)^\prime = 0,
\end{align}
\end{subequations} 
for $i = 1,2$ with domain $[-1,1]$.
The coefficients are
\begin{equation}
\chi_1 = \frac{e \phi_{0}}{k_{B} T}, \quad 
\chi_ 2 = \frac{e c_{0} L^{2}}{\phi_{0} \epsilon_{t}},
\end{equation}
where 
$\chi_{1}$ is the ratio of electric potential energy 
to thermal energy,
and
$\chi_{2}$ is the inverse square of the scaled Debye length.
Except as otherwise noted,
the non-dimensional diffusion coefficient is $D_i=1$.


\subsection{Problem specification}

The first set of examples are
electroneutral and non-electroneutral systems,
where the 
non-dimensional PNP equations~\eqref{eqn:PNP_system} 
are solved on the interval $[-1,1]$.
We consider Robin boundary conditions for the potential,
\begin{equation}
\label{robin4P}
\phi(x)\pm\eta \phi^\prime(x) =\phi_{\pm}, \quad x =-1,1, 
\end{equation}
where $\eta, \phi_{\pm}$ are constants that depend on
biological modeling considerations.
Note that if $\eta=0$, 
then~\eqref{robin4P} reduces to Dirichlet boundary conditions
for the potential.
The channel is assumed to be closed at both ends,
hence the no-flux boundary condition is applied
for the ion concentrations,
\begin{equation}
\label{nofluxBC4NP}
D_i\left(c_i^\prime(x) + \chi_1 z_i c_i(x)\phi^\prime(x) \right) = 0, 
\quad x = -1,1.
\end{equation}
Finally, the total concentration for each ion species 
is a specified positive constant,
\begin{equation}
\label{normalize4NP}
\int_{-1}^{1}c_i(x)dx = a_i, \quad i = 1,2.
\end{equation}
In addition to the 
electroneutral and non-electroneautral systems,
we also consider a K$^+$~ion channel model
and
the problem specification in that case
is given further below.


\section{Integral form}

This section recasts the 
1D steady state PNP differential equations 
as a set of integral equations in the form
\begin{subequations}
\label{eqn:map}
\begin{align}
\label{eqn:map_P}
\phi^\prime &= {\bf P}[c_1, c_2], \\[5pt]
\label{eqn:map_NP}
c_i &= {\bf NP}[c_i, \phi^\prime], \quad i = 1,2,
\end{align}
\end{subequations}
where ${\bf P}, {\bf NP}$ are operators
obtained by integrating the Poisson 
and Nernst-Planck equations
while incorporating the Robin 
and 
no-flux boundary conditions
and
the total concentration condition.
Note that~\eqref{eqn:map} defines a fixed-point problem,
where ${\bf P}$ maps
the concentrations $c_1, c_2$ into 
the potential gradient $\phi^\prime$,
while ${\bf NP}$ maps the concentration $c_i$ 
and 
potential gradient $\phi^\prime$ into
the concentration $c_i$ for $i=1,2$.
The integral form of the 
PNP equations in~\eqref{eqn:map}
is solved by Gummel iteration~\cite{gummel1964self}
for $\phi^\prime, c_1, c_2$,
and
once it converges,
further processing yields 
the potential~$\phi$
and
the concentration gradients $c_1^\prime, c_2^\prime$
needed to compute the current.
Note that a fixed-point formulation
and 
Gummel iteration 
was also used in the analysis
of semiconductor drift-diffusion equations,
although without recasting the differential equations as 
integral equations~\cite{jerome1985consistency}.

The derivation of the 
${\bf P}$ and ${\bf NP}$
operators relies on Green's 3rd~identity,
\begin{equation}\label{2ndIdentity}
\begin{aligned}
f(x) 
&= \Big[g(x,y)f^\prime(y) - g_y(x,y)f(y)\Big]_{y=-1}^{y=1} \\
&- \int_{-1}^1 g(x,y)f^{\prime\prime}(y) dy,
\end{aligned}
\end{equation}
where $f(x)$ is twice continuously differentiable on $[-1,1]$
and
$g(x,y)$ is the 
1D free-space Laplace Green's function,
\begin{equation}\label{green1D}
g(x,y) = 
-\frac{1}{2}\lvert x-y\rvert.
\end{equation}
Note that among the two terms 
on the right side of~\eqref{2ndIdentity},
the first term (or bracket term)
involves the boundary values 
$f(\pm 1), f^\prime(\pm 1)$,
and
the second term 
(or volume integral term)
involves $f^{\prime\prime}(x)$ over the entire interval.
One can think of~\eqref{2ndIdentity}
as defining an operator,
where $f(\pm 1),f^\prime(\pm 1),f^{\prime\prime}(x)$
are the input
and
$f(x)$ is the output.
This point of view is taken in reformulating
the PNP equations.


\subsection{Derivation of {\bf P} operator}
\label{integral_form_of_P_equation}

Here we assume the ion concentrations
$c_1(x), c_2(x)$ are given
and
we shall derive an integral expression for the 
potential gradient 
$\phi^\prime(x)$.
The Poisson equation~\eqref{eqn:P_non_dimensional} 
is written as
\begin{equation}
\label{eqn:P}
\phi^{\prime\prime}(x) = 
-\frac{\chi_2}{\varepsilon}
\sum\limits_{i=1}^2 z_i c_i(x),
\end{equation}
and then applying~\eqref{2ndIdentity}
with $f(y) = \phi(y)$ yields an 
integral expression for the potential,
\begin{equation}
\label{integral_P}
\begin{aligned}
\phi(x) & = 
\Big[g(x,y) \phi^\prime(y) -g_y(x,y)\phi(y)\Big]_{y=-1}^{y=1} \\
& + \frac{\chi_2}{\epsilon}
\sum_{i=1}^2 z_i\int_{-1}^1 g(x,y)c_i(y)dy.
\end{aligned}
\end{equation}
Differentiating~\eqref{integral_P}
yields an analogous expression for the potential gradient,
\begin{equation}
\label{integral_P_gradient}
\begin{aligned}
\phi^\prime(x) 
& = \Big[g_x(x,y) \phi^\prime(y) -
g_{xy}(x,y)\phi(y)\Big]_{y=-1}^{y=1} \\
& +\frac{\chi_2}{\epsilon}
\sum_{i=1}^2 z_i\int_{-1}^1 g_x(x,y)c_i(y)dy.
\end{aligned}
\end{equation}
The bracket terms in~\eqref{integral_P}
and
\eqref{integral_P_gradient}
involve
the boundary values of the potential $\phi(\pm 1)$
and
potential gradient $\phi^\prime(\pm 1)$.
The Robin boundary conditions~\eqref{robin4P}
provide two equations for these four unknowns
and
two more equations are obtained as follows.

Integrating the 
Poisson equation \eqref{eqn:P}
over $[-1,x]$ yields
\begin{equation}
\label{eqn:P_intermediate}
\phi^\prime(x) - \phi^\prime(-1) =
-\frac{\chi_2}{\epsilon}
\int_{-1}^x
\sum_{i=1}^2 z_i c_i(s) ds,
\end{equation}
and
setting $x=1$ yields
\begin{equation}
\label{eqn:P_BC1}
\phi^\prime(1) - \phi^\prime(-1) =
-\frac{\chi_2}{\epsilon}
\sum_{i=1}^2 z_i a_i.
\end{equation}
Then integrating \eqref{eqn:P_intermediate}
over $[-1,1]$ yields
\begin{equation}
\label{eqn:P_BC_2a}
\phi(1) - \phi(-1) - 2\phi^\prime(-1) =
-\frac{\chi_2}{\epsilon}
\sum_{i=1}^2 z_i (a_i - b_i),
\end{equation}
where $a_i$ is the total concentration of ion species~$i$
defined in~\eqref{normalize4NP} and
\begin{equation}
\label{eqn:concentration_first_moment}
b_i = \int_{-1}^1 xc_i(x) dx, \quad i = 1,2,
\end{equation}
is the first moment of the ion concentration.
By similar steps,
integrating the Poisson equation \eqref{eqn:P} 
over $[x,1]$ yields
\begin{equation}
\label{eqn:P_BC_2b}
\phi(1) - \phi(-1) -2\phi^\prime(1) = 
\frac{\chi_2}{\epsilon}\sum_{i=1}^2 z_i (a_i + b_i).
\end{equation}
Then adding 
\eqref{eqn:P_BC_2a} and \eqref{eqn:P_BC_2b} yields
\begin{equation}
\label{eqn:P_BC2}
\phi(1) - \phi(-1) 
-\phi^\prime(1) - \phi^\prime(-1) =
\frac{\chi_2}{\epsilon}\sum_{i=1}^2 z_i b_i.
\end{equation}
Combining the 
Robin boundary conditions~\eqref{robin4P}
with~\eqref{eqn:P_BC1} and~\eqref{eqn:P_BC2}
yields the linear system
\begin{equation}
\label{eqn:P_system}
\begin{pmatrix}
1 & \eta & 0 & 0 \\
0 & 0 & 1 & -\eta \\
0 & 1 & 0 & -1 \\
1 & -1 & -1 & -1 \\
\end{pmatrix}
\begin{pmatrix}
\phi(1) \\
\phi^\prime(1) \\
\phi(-1) \\
\phi^\prime(-1)
\end{pmatrix} =
\begin{pmatrix}
\phi_+ \\
\phi_- \\
\displaystyle -\frac{\chi_2}{\epsilon}
\sum\limits_{i=1}^2 z_i a_i \\
\displaystyle \hfill \frac{\chi_2}{\epsilon}
\sum\limits_{i=1}^2 z_i b_i \\
\end{pmatrix}.
\end{equation}
Hence given ion concentrations 
$c_1(x), c_2(x)$ for $x \in [-1,1]$,
the ${\bf P}$ operator in~\eqref{eqn:map_P}
is defined by the following steps.

\begin{itemize}
\item
compute the side right of the linear system~\eqref{eqn:P_system}
and
solve for the potential boundary values 
$\phi(\pm 1), \phi^\prime(\pm 1)$
\item
compute the interior values of the
potential gradient $\phi^\prime(x)$ for $x \in (-1,1)$
from~\eqref{integral_P_gradient}
\end{itemize} 
These steps require computing certain integrals
and
the discretization schemes will be explained further below.


\subsection {Derivation of {\bf NP} operator} 
\label{integral_form_of_NP_equation}

Here we assume the ion concentration $c_i(x)$
and
potential gradient $\phi^\prime(x)$ are given
and
we shall derive an integral expression for $c_i(x)$.
The Nernst-Planck 
equation~\eqref{eqn:NP_non_dimensional} is written as
\begin{equation}
c_i^{\prime\prime}(x) = -\chi_1 z_i (c_i(x)\phi^\prime(x))^\prime,
\end{equation}
after cancelling the constant diffusion coefficient $D_i$,
and
then applying~\eqref{2ndIdentity} 
with $f(x) = c_i(x)$ yields an
integral expression for the ion concentration,
\begin{equation}\label{integral_NP}
\begin{aligned}
c_i(x) &= \Big[g(x,y) c_i^\prime(y) -
g_y(x,y)c_i(y)\Big]_{y=-1}^{y=1} \\ 
&+ \chi_1z_i\int_{-1}^{1} g(x,y) \left(c_i(y)\phi^\prime(y)\right)^\prime dy.
\end{aligned}
\end{equation}
Differentiating~\eqref{integral_NP}
yields an analogous expression for the concentration gradient,
\begin{equation}\label{integral_NP_gradient}
\begin{aligned}
c_i^\prime(x) &= \Big[g_x(x,y) c_i^\prime(y) -
g_{xy}(x,y)c_i(x)\Big]_{y=-1}^{y=1} \\ 
&+ \chi_1z_i\int_{-1}^{1} g_x(x,y)  \left(c_i(y)\phi^\prime(y)\right)^\prime dy.
\end{aligned}
\end{equation}
The bracket terms in~\eqref{integral_NP} 
and \eqref{integral_NP_gradient}
involve the boundary values of the
concentration $c_i(\pm 1)$ and
concentration gradient $c_i^\prime(\pm 1)$.
The solution procedure in this case is 
somewhat different than
what was described above for the ${\bf P}$ operator.

From the Nernst-Planck equation~\eqref{eqn:NP_non_dimensional}
and
no-flux boundary conditions~\eqref{nofluxBC4NP},
it follows that
\begin{equation}
\label{eqn:NP_new}
c_i^\prime(x) + \chi_1 z_i c_i(x)\phi^\prime(x) = 0, \quad x \in (-1,1).
\end{equation}
Integrating~\eqref{eqn:NP_new} over $[-1,x]$ yields
\begin{equation}
\label{eqn:NP}
c_i(x) - c_i(-1) =
-\chi_1 z_i \int_{-1}^x c_i(s)\phi^\prime(s)ds,
\end{equation}
and
setting $x = 1$ yields
\begin{equation}
\label{eqn:NP_1}
c_i(1) - c_i(-1) =
-\chi_1 z_i \int_{-1}^1 c_i(s)\phi^\prime(s)ds.
\end{equation}
Then integrating~\eqref{eqn:NP} over $[-1,1]$ yields
\begin{equation}
\label{eqn:NP_2}
c_i(1) + c_i(-1) =
a_i - \chi_1z_i\!\int_{-1}^1 \!\! xc_i(x)\phi^\prime(x)dx.
\end{equation}
Hence given the ion concentration $c_i(x)$ 
and
potential gradient $\phi^\prime(x)$ for $x \in [-1,1]$,
the ${\bf NP}$ operator in~\eqref{eqn:map_NP} is defined
by the following steps.
\begin{itemize}
\item
compute the integrals on the 
right of~\eqref{eqn:NP_1}-\eqref{eqn:NP_2}
and
solve for the concentration boundary values $c_i(\pm 1)$
on the left
\item
compute the concentration gradient 
boundary values $c_i^\prime(\pm 1)$ using $c_i(\pm 1)$
and
the no-flux boundary conditions~\eqref{nofluxBC4NP}
\item
compute the interior concentration $c_i(x)$ for $x \in (-1,1)$
from~\eqref{integral_NP}
\end{itemize}
This concludes the derivation of the
${\bf P}, {\bf NP}$ operators
in the integral form of the 
PNP equations~\eqref{eqn:map}.


\section {Numerical method} 

This section presents the numerical method used to
solve the integral form of the PNP equations
with the 
Robin and no-flux boundary conditions
and
total concentration constraint.
The domain $[-1,1]$ is divided into $N$ subintervals 
with endpoints $x_k, k=0:N$, 
where $x_0 = -1, x_N = 1$,
and
the grid spacing is $h_k = x_{k+1}-x_k$.
We consider two point sets, 
uniform and Chebyshev,
\begin{subequations}
\begin{align}
x_k & = -1 + 2k/N, \quad k = 0:N, \\[5pt]
x_k & = -\cos\theta_k, \quad \theta_k = k\pi/N, \quad k = 0:N.
\end{align}
\end{subequations}
The potential and concentration values 
and
their gradients are denoted as
\begin{subequations}
\label{variables}
\begin{align}
\phi_k &= \phi(x_k), \quad 
\phi^\prime_k = \phi^\prime(x_k) \\
c_{i,k} &= c_i(x_k), \quad 
c_{i,k}^\prime = c_i^\prime(x_k).
\end{align}
\end{subequations}

Two types of integrals 
appear in the ${\bf P}, {\bf NP}$ operators 
depending on whether or not the Green's function is involved.
When the Green's function is not involved
the integral is discretized by the trapezoid rule,
\begin{subequations}
\label{eqn:total_concentration_discrete}
\begin{align}
b_i &= 
\int_{-1}^1 xc_i(x) dx \approx \sum_{k=0}^N x_k c_{i,k}w_k, \\
&\int_{-1}^1 c_i(x)\phi^\prime(x)dx\approx \sum_{k=0}^N  c_{i,k}\phi^\prime_k w_k, \\
&\int_{-1}^1 xc_i(x)\phi^\prime(x)dx\approx \sum_{k=0}^N x_k c_{i,k}\phi^\prime_kw_k,
\end{align}
\end{subequations}
where $w_k$ are the trapezoid weights.
Next we explain how the integrals involving the 
Green's function are computed.


\subsection{Integrals for {\bf P} operator}
\label{section:discreteP}

First consider the integral in
the potential expression~\eqref{integral_P}, 
\begin{subequations}
\begin{align}
\int_{-1}^1 &g(x,y)c_i(y)dy \\ 
&=
\sum_{j=0}^{N-1}\int_{x_j}^{x_{j+1}}g(x,y)c_i(y)dy \\
\label{integral_P_midpoint}
&\approx
\sum_{j=0}^{N-1}g(x,x_{j+1/2})\int_{x_j}^{x_{j+1}}c_i(y)dy \\
\label{integral_P_trapezoid}
&\approx
\sum_{j=0}^{N-1}g(x,x_{j+1/2})
\frac{c_{i,j}+c_{i,j+1}}{2}h_j,
\end{align}
\end{subequations}
where $x_{j+1/2} = (x_j+x_{j+1})/2$.
Note that a midpoint approximation 
is used in~\eqref{integral_P_midpoint}
and
the trapezoid rule is used in~\eqref{integral_P_trapezoid}.
Then setting $x = x_k$ in~\eqref{integral_P} yields
\begin{equation}
\label{eqn:potential_discrete}
\begin{aligned}
& \phi_k = \Big[g(x_k, y) \phi^\prime(y) -
g_y(x_k,y)\phi(y)\Big]_{y=-1}^{y=1} \\
& +\frac{\chi_2}{\epsilon}
\sum\limits_{i=1}^2 z_i
\sum_{j=0}^{N-1}
g(x_k,x_{j+1/2})\frac{c_{i,j}+c_{i,j+1}}{2}h_j.
\end{aligned}
\end{equation}
Applying the same steps to the
integral in the 
potential gradient exression~\eqref{integral_P_gradient} 
yields
\begin{equation}
\label{eqn:potential_gradient_discrete}
\begin{aligned}
& \phi_k^\prime = \Big[g_x(x_k, y) \phi^\prime(y) -
g_{xy}(x_k, y)\phi(y)\Big]_{y=-1}^{y=1} \\
& +\frac{\chi_2}{\epsilon}
\sum\limits_{i=1}^2 z_i
\sum_{j=0}^{N-1}
g_x(x_k, x_{j+1/2})\frac{c_{i,j}+c_{i,j+1}}{2}h_j.
\end{aligned}
\end{equation}
Equations~\eqref{eqn:potential_discrete}
and~\eqref{eqn:potential_gradient_discrete}
are used to compute the 
interior values of the potential and its gradient,
$\phi_k, \phi_k^\prime, k=1:N-1$,
where the boundary values 
$\phi(\pm 1), \phi^\prime(\pm 1)$
in the bracket terms
are computed as described in 
Section~\ref{integral_form_of_P_equation}.


\subsection{Integrals for {\bf NP} operator}
\label{section:discreteC}

Now consider the integral in
the concentration expression~\eqref{integral_NP},
\begin{subequations}
\label{integral_NP_discrete}
\begin{align}
&\int_{-1}^{1} g(x,y) 
\left(c_i(y)\phi^\prime(y)\right)^\prime dy \\
&= 
\sum_{j=0}^{N-1}\int_{x_j}^{x_{j+1}}
g(x, y)\left(c_i(y)\phi^\prime(y)\right)^\prime dy \\
\label{integral_NP_midpoint}
&\approx  
\sum_{j=0}^{N-1}g(x,x_{j+1/2})\int_{x_j}^{x_{j+1}}
\left(c_i(y)\phi^\prime(y)\right)^\prime dy \\
\label{integral_NP_FTC}
&=
\sum_{j=0}^{N-1} g(x, x_{j+1/2})
(c_{i,j+1}\phi^\prime_{j+1} - c_{i,j}\phi^\prime_j),
\end{align}
\end{subequations}
where a midpoint approximation is used 
in~\eqref{integral_NP_midpoint}
and
the fundamental theorem of calculus is used in~\eqref{integral_NP_FTC}.
Then setting $x = x_k$ 
in~\eqref{integral_NP} yields
\begin{equation}
\label{eqn:concentration_discrete}
\begin{aligned}
& c_{i,k} = \Big[g(x_k,y) c_i^\prime(y) - g_y(x_k, y)c_i(y)\Big]_{y=-1}^{y=1} \\ 
& + \chi_1z_i\sum_{j=0}^{N-1} g(x_k,x_{j+1/2})
(c_{i,j+1}\phi^\prime_{j+1} - c_{i,j}\phi^\prime_j).
\end{aligned}
\end{equation}
Applying the same steps to the integral in the
concentration gradient 
expression~\eqref{integral_NP_gradient} yields
\begin{equation}
\label{eqn:concentration_gradient_discrete}
\begin{aligned}
& c_{i,k}^\prime = 
\Big[g_x(x_k, y) c_i^\prime(y)-
g_{xy}(x_k, y)c_i(y)\Big]_{y=-1}^{y=1} \\ 
& + \chi_1z_i\sum_{j=0}^{N-1} g_x(x_k,x_{j+1/2})
(c_{i,j+1}\phi^\prime_{j+1} - c_{i,j}\phi^\prime_j).
\end{aligned}
\end{equation}
Equations~\eqref{eqn:concentration_discrete}
and~\eqref{eqn:concentration_gradient_discrete}
are used to compute the 
interior values of the 
concentration and its gradient,
$c_{i,k}, c_{i,k}^\prime, k=1:N-1$,
where the boundary values 
$c_i(\pm 1), c_i^\prime(\pm 1)$
in the bracket terms
are computed as described in 
Section~\ref{integral_form_of_NP_equation}.


\subsection{Gummel iteration}

Define vectors containing the 
grid values of the potential, ion concentration,
and their gradients,
\begin{subequations}
\begin{align}
{\bf \Phi} 
&= [\phi_0, \ldots, \phi_N]^T, \\
{\bf D\Phi} 
&= [\phi_0^\prime, \ldots, \phi_N^\prime]^T, \\
{\bf c}_i 
&= [c_{i,0}, \ldots, c_{i,N}]^T, \quad i=1,2, \\
{\bf Dc}_i 
&= [c_{i,0}^\prime, \ldots, c_{i,N}^\prime]^T, \quad i=1,2.
\end{align}
\end{subequations}
Then after discretizing,
the integral form of the PNP equations~\eqref{eqn:map} 
becomes a system of nonlinear equations for the grid values,
\begin{subequations}
\label{eqn:discrete}
\begin{align}
\label{eqn:map_P_discrete}
{\bf D\Phi} &= {\bf P}[{\bf c}_1, {\bf c}_2], \\[5pt]
\label{eqn:map_NP_discrete}
{\bf c}_i &= {\bf NP}[{\bf c}_i, {\bf D\Phi}], \quad i=1,2.
\end{align}
\end{subequations}
To solve~\eqref{eqn:discrete}
we employ Gummel iteration
with relaxation~\cite{gummel1964self}.
Letting superscript $n$ denote the iteration step,
we have
\begin{subequations}
\label{eqn:Gummel}
\begin{align}
\label{eqn:Gummel_a}
{\bf D\Phi}^*
&= {\bf P}[{\bf c}_1^{(n)}, {\bf c}_2^{(n)}], \\[2.5pt]
\label{eqn:Gummel_b}
{\bf D\Phi}^{(n+1)} 
&= \omega {\bf D\Phi}^* + 
(1-\omega) {\bf D\Phi}^{(n)}, \\[2.5pt]
\label{eqn:Gummel_c}
{\bf c}_i^* 
&= {\bf NP}[{\bf c}_i^{(n)}, {\bf D\Phi}^{(n+1)}],
\quad i = 1,2, \\[2.5pt]
\label{eqn:Gummel_d}
{\bf c}_i^{(n+1)} 
&= \omega {\bf c}_i^* + 
(1-\omega) {\bf c}_i^{(n)}, \quad i=1,2,
\end{align}
\end{subequations}
where the relaxation parameter is chosen with
$0 < \omega < 1$.
The intermediate potential gradient 
${\bf D\Phi}^*$ 
is computed 
by the procedure in Section~\ref{section:discreteP},
and
the intermediate concentration 
${\bf c}_i^*, i=1,2$ 
is computed
by the procedure in Section~\ref{section:discreteC}.
Hence each iteration has two stages;
the first stage 
comprising~\eqref{eqn:Gummel_a}-\eqref{eqn:Gummel_b}
computes the 
potential gradient using the 
ion concentrations from the previous iteration,
and
the second stage 
comprising~\eqref{eqn:Gummel_c}-\eqref{eqn:Gummel_d}
computes the 
ion concentrations using the
potential gradient from the first stage.

The scheme starts with constant initial guesses 
for the potential gradient 
and 
ion concentrations,
\begin{subequations}
\begin{align}
\label{eqn:initial_phi}
& ({\bf D\Phi}^{(0)})_j =
\frac{\phi_+ - \phi_-}{2}, \\ 
\label{eqn:initial_c}
& ({\bf c}_i^{(0)})_j = \frac{1}{2}a_i, \quad i=1,2,
\end{align}
\end{subequations}
where $j = 0:N$ is the grid-point index,
the right side of
\eqref{eqn:initial_phi} 
is the gradient determined by the 
formal potential boundary values $\phi_{\pm}$,
and
\eqref{eqn:initial_c} 
satisfies the total concentration constraint~\eqref{normalize4NP}.
The stopping criterion is
\begin{subequations}
\label{eqn:stopping_criterion}
\begin{align}
\lvert\lvert{\bf D\Phi}^{(n+1)}-{\bf D\Phi}^{(n)}\rvert\rvert_2 < 10^{-6}, \\
\max\limits_{i=1,2}\, 
\lvert\lvert{\bf c}_i^{(n+1)}-{\bf c}_i^{(n)}\rvert\rvert_2 < 10^{-6}.
\end{align}
\end{subequations}
After the scheme converges,
the grid values of the potential $\phi_k$
are computed from~\eqref{eqn:potential_discrete},
and
the grid values of the 
concentration gradient $c_{i,k}^\prime$ 
are computed 
from~\eqref{eqn:concentration_gradient_discrete}.
Note that the potential is presented
for diagnostic purposes
and
the concentration gradient is
used to compute the current
in the K$^+$ channel simulation further below.

This completes the description of the 
proposed integral equation method for the
1D steady state PNP equations.
Several antecedents have used integral equation
formulations to solve 
two-point boundary value problems
although to the best of our knowledge,
they used Green's functions 
that satisfy homogeneous boundary conditions~\cite{
barcilon1992ion,
greengard1991numerical,
viswanath2013navier},
whereas the present PNP solver uses the 
free-space Green's function which we expect 
will better facilitate a future extension to 2D and 3D problems.


\section{Numerical results}

In this section the integral method is applied to compute
several examples from prior work 
in which the PNP equations were solved by
finite-difference or finite-element methods~\cite{flavell2014conservative,lee2010new}.
In these examples there are two ion species 
with 
valence $z_1=-1$ (anions), $z_2=1$ (cations),
and
two types of systems
depending on the total ion concentration,
electroneutral ($a_1 = a_2$)
and
non-electroneutral ($a_1 \ne a_2$).
Table~\ref{tab:parameters} 
gives the parameter values.
The calculations used either 
Chebyshev or uniform grid points as noted
and
the relaxation parameter $\omega$ 
was chosen by trial and error.
The method was programmed in Python
and
the calculations were done on an iMac computer 
with the Apple M1 chip.

\begin{table*}[t]
\centering
\begin{tabular}{|c|c|c|c|c|c|c|c|c|c|c|}
\hline
case & ref & $a_1$ & $a_2$  & $\phi_-$ & $\phi_+$ & $\epsilon$ 
& $\eta$ & $\chi_1$ & $\chi_2$ & $\omega$ \\ 
\hline
\multicolumn{11}{|c|}{electroneutral systems} \\
\hline
1.1 & \cite{flavell2014conservative},\cite{lee2010new} & 1 & 1  & -1 &  1 & 1/4    & $\epsilon$   & 1   &     $1/\epsilon$ & 0.7  \\ \hline
1.2 & \cite{flavell2014conservative},\cite{lee2010new} & 1 & 1  & -1 &  1 & 1/64   & $\epsilon$   & 1   &    $1/\epsilon$  & 0.09 \\ \hline
2.1 & \cite{lee2010new} & 1 & 1  & -1 &  1 & 1/4    & 
$0, \epsilon^2, \epsilon, \epsilon^{1/2}, 1$  & 1   &     $1/\epsilon$ &  0.7 \\  \hline
2.2 & \cite{lee2010new} & 1 & 1  & -1 &  1 & 1/64  & 
$0, \epsilon^2, \epsilon, \epsilon^{1/2}, 1$ & 1   &    $1/\epsilon$  & 0.09 \\  \hline
3 & \cite{flavell2014conservative} & 2 & 2   &  1 & -1  & 1     & 4.63e-5 & 3.1 & 125.4 & 0.6 \\ \hline
\multicolumn{11}{|c|}{non-electroneutral systems} \\ \hline
4.1 & \cite{lee2010new} & 1 & 2      & 1 & 1  & 1/4 & 
$0, \epsilon^2, \epsilon, \epsilon^{1/2}, 1$  & 1   & $1/\epsilon$ & 0.6 \\  \hline
4.2 & \cite{lee2010new} & 1 & 2   & 1 & 1  & 1/16  & 
$0, \epsilon^2, \epsilon, \epsilon^{1/2}, 1$ & 1   & $1/\epsilon$  & 0.16 \\  \hline
\end{tabular}
\vskip 5pt
\caption{
Parameters for electroneutral and
non-electroneutral systems,
case number,
reference,
total concentration of anions and cations $a_1, a_2$,
Robin boundary condition potential parameters $\phi_\pm$,
permittivity $\epsilon$,
capacitance parameter $\eta$,
energy ratio $\chi_1$,
inverse square Debye length $\chi_2$,
relaxation parameter $\omega$.}
\label{tab:parameters}
\end{table*}


\subsection{Electroneutral systems}

Here we consider three electroneutral systems
in which the 
total anion and cation concentrations are equal, $a_1 = a_2$,
the formal potential boundary values are
either $\phi_- = -1, \phi_+ = 1$ (case 1, 2)
or
$\phi_- = 1, \phi_+ = -1$ (case 3).


\subsubsection{Electroneutral case 1}

\begin{figure*}[htb!]
\centering
\includegraphics[width=\textwidth]{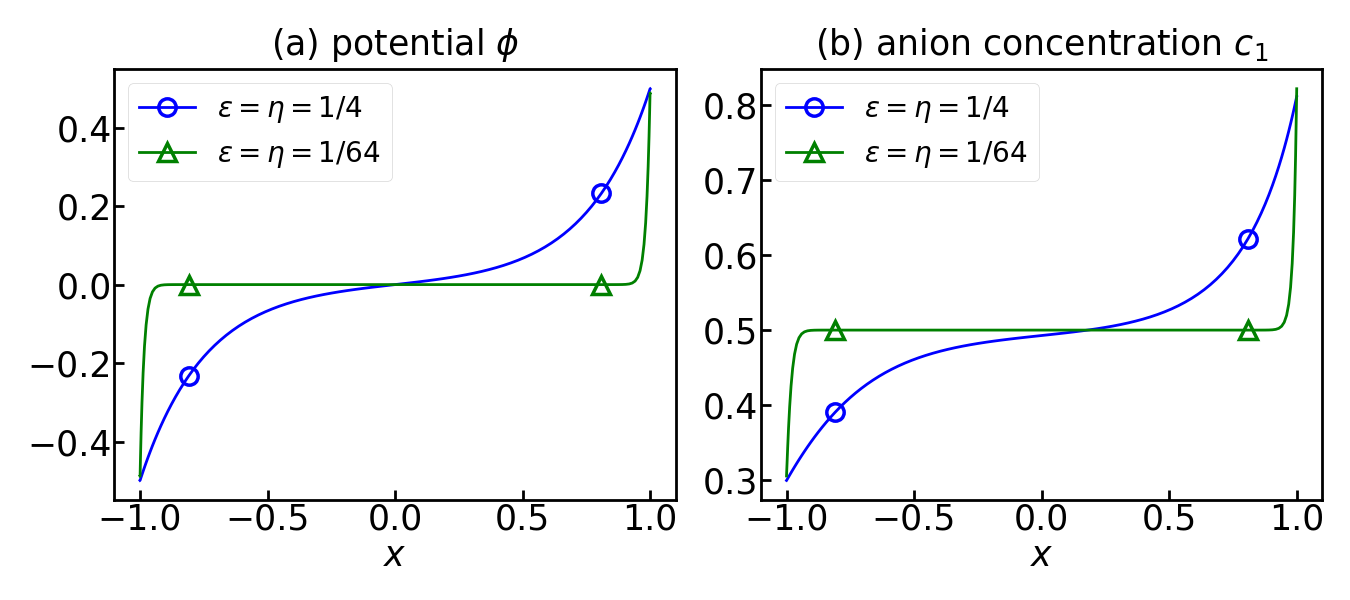}
\caption{
Electroneutral systems,
case 1.1, 1.2 in Table~\ref{tab:parameters},
large permittivity $\epsilon = 1/4$,
small permittivity $\epsilon = 1/64$,
(a) potential $\phi$,
(b) anion concentration $c_1$,
integral method with Chebyshev points,
$N=100$ subintervals.}
\label{fig:CompareCaseNeutral1}
\end{figure*}

In this case
the capacitance parameter is equal to the permittivity,
$\eta = \epsilon$,
and
two values are considered,
$\epsilon = 1/4, 1/64$~\cite{flavell2014conservative,lee2010new}.
Figure~\ref{fig:CompareCaseNeutral1} 
shows the computed potential $\phi$
and 
anion concentration $c_1$;
the cation concentration $c_2$ is not shown,
since it is simply the reflection of $c_1$
about the channel midpoint,
$c_2(x) = c_1(-x)$.
With $\phi_\pm = \pm 1$,
the potential is an odd function,
$\phi(-x) = -\phi(x)$,
and
it increases monotonically from left to right 
across the channel.
The anions are repelled from the low potential at $x=-1$
and
are attracted to the high potential at $x=1$
where they accumulate due to the no-flux boundary condition.
For large permittivity $\epsilon = 1/4$
the potential and concentration profiles have mild gradients,
but for small permittivity $\epsilon = 1/64$
the profiles have sharp boundary layers.
These results agree well with 
Fig.~1 in~\cite{flavell2014conservative,lee2010new}.

Following~\cite{flavell2014conservative}
the convergence rate for the potential is defined by
\begin{equation}
p = -\log_2 
\frac{\|{\bf \Phi}_{4N}-{\bf \Phi}_{2N}\|_2}
{\|{\bf \Phi}_{2N}-{\bf \Phi}_{N}\|_2},
\end{equation}
where ${\bf \Phi}_N$ denotes the vector of
computed potential grid values using $N$ subintervals.
Table~\ref{tab:Electroneutral1}
shows results for both uniform and Chebyshev points
indicating that $p \approx 2$,
so the integral method is 2nd order accurate.
For large permittivity $\epsilon = 1/4$ where the solution is smooth,
both point sets require the same number of iterations,
while the Chebyshev points yield slightly smaller error.
For small permittivity $\epsilon = 1/64$
where the solution has a sharp boundary layer,
the number of iterations is higher for both point sets,
but the Chebyshev points require fewer iterations
and
yield smaller error than the uniform points.


\begin{table}
\centering
\begin{tabular}{|c|c|c|c|c|c|c|} 
\hline
\multicolumn{7}{|l|}{case 1.1, $\epsilon = 1/4$} \\
\hline
\multicolumn{1}{|c}{$N$} & \multicolumn{3}{|c|}{uniform} & \multicolumn{3}{c|}{Chebyshev} \\
\hline
& iter & error & $p$ & iter & error & $p$ \\
\hline
 50 & 19 & 2.1e-3 &       & 19 & 2.1e-3 & \\ \hline
100 & 17 & 5.3e-4 & 2.007 & 17 & 4.3e-4 & 1.997 \\ \hline
200 & 17 & 1.3e-4 & 2.002 & 17 & 3.8e-5 & 2.003 \\ \hline
400 & 17 & 3.3e-5 & 2.000 & 17 & 1.8e-5 & 2.000 \\ \hline
\hline
\multicolumn{7}{|l|}{case 1.2, $\epsilon = 1/64$} \\ 
\hline
\multicolumn{1}{|c}{$N$} & \multicolumn{3}{|c|}{uniform} & \multicolumn{3}{c|}{Chebyshev} \\
\hline
& iter & error & $p$ & iter & error & $p$ \\
\hline
 50 & nc    &        &       & 245 & 5.1e-2 &       \\ \hline
100 & 21008 & 2.2e-1 &       & 202 & 1.0e-2 & 1.871 \\ \hline
200 & 1180  & 2.8e-2 & 2.452 & 199 & 2.8e-3 & 2.027 \\ \hline
400 & 656   & 6.3e-3 & 2.098 & 206 & 7.4e-4 & 2.003 \\ \hline
\end{tabular}
\vskip 5pt
\caption{Electroneutral system
(case 1.1, 1.2),
parameters in Table~\ref{tab:parameters},
large permittivity $\epsilon = 1/4$,
small permittivity $\epsilon = 1/64$,
$N$ subintervals, uniform and Chebyshev points,
number of iterations iter, 
potential error $\|{\bf \Phi}_{2N}-{\bf \Phi}_{N}\|_2$, 
convergence rate $p$,
nc = did not converge.
}
\label{tab:Electroneutral1}
\end{table}


\subsubsection{Electroneutral case 2}

\begin{figure*}[htb!]
\centering
\includegraphics[width=\textwidth]{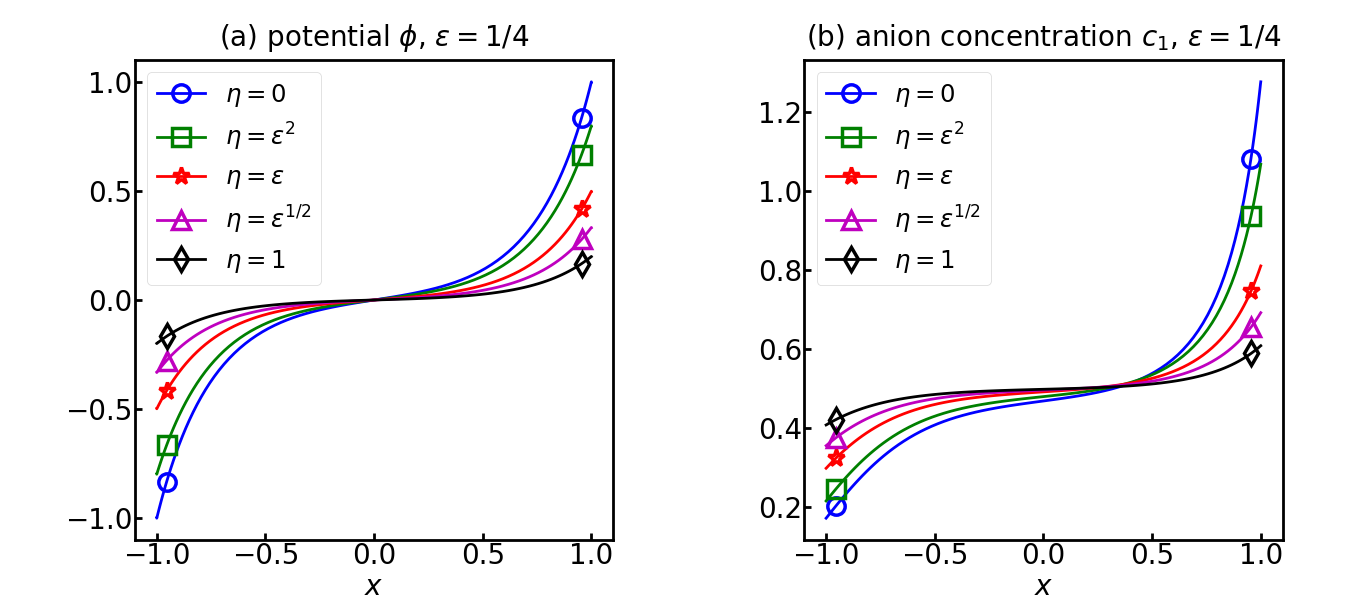}
\includegraphics[width=\textwidth]{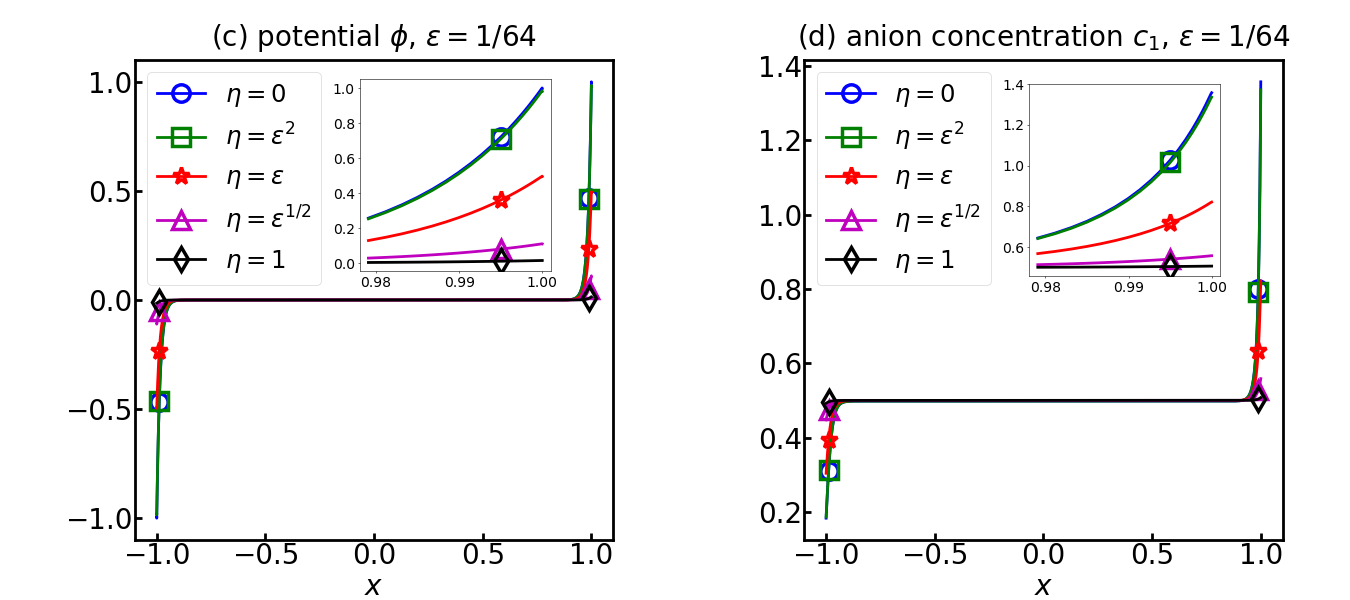}
\caption{
Electroneutral systems, 
case 2.1, 2.2 in Table~\ref{tab:parameters},
potential $\phi$,
anion concentration $c_1$,
large permittivity $\epsilon = 1/4$ (a,b),
small permitivity $\epsilon = 1/64$ (c,d),
integral method with Chebyshev points,
$N=100$ subintervals.}
\label{fig:CaseCompare23}
\end{figure*}

This case extends case 1 by considering 
five values of the capacitance parameter,
$\eta = 0, \epsilon^2, \epsilon, \epsilon^{1/2}, 1$.
Figure~\ref{fig:CaseCompare23} 
shows the potential $\phi$ 
and 
anion concentration $c_1$,
where the computed boundary values of $\phi(1), c_1(1)$ 
are printed alongside the curves
for comparison with~\cite{lee2010new}.
As in case 1,
the anions are repelled from the 
potential minimum at $x=-1$
and
are attracted to the potential maximum at $x=1$;
moreover
for large permittivity $\epsilon = 1/4$
the profiles have mild gradients,
but for small permittivity $\epsilon = 1/64$
they have sharp boundary layers.
The dependence on $\eta$ is interesting.
As $\eta \to 0$,
the potential boundary values converge
to the formal values, $\phi(\pm 1) \to \phi_\pm$,
as required by the 
Robin boundary conditions~\eqref{robin4P}.
On the other hand as $\eta$ increases,
the potential vanishes, $\phi(x) \to 0$,
and
the anion concentration tends to the 
bulk value, $c_1(x) \to 1/2$;
these features are consistent with the
Robin and no-flux 
boundary conditions~\eqref{robin4P},\eqref{nofluxBC4NP}
and
the total concentration condition~\eqref{normalize4NP}.
The present results agree well with 
Fig.~1 in~\cite{lee2010new}.

\subsubsection{Electroneutral case 3}

\begin{figure*}[t]
\centering
\includegraphics[width=\textwidth]{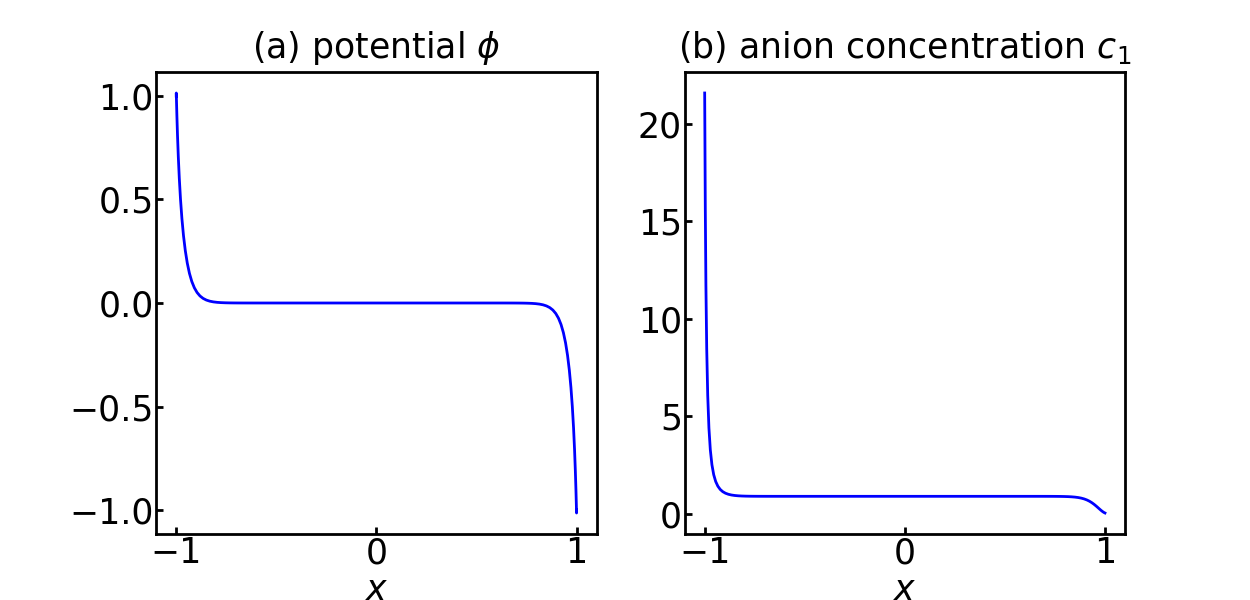}
\caption{
Electroneutral system, 
case 3 in Table~\ref{tab:parameters},
(a) potential $\phi$,
(b) anion concentration $c_1$,
integral method with Chebyshev points,
$N=100$ subintervals.}
\label{fig:CompareCaseNeutral2}
\end{figure*}

This case uses the more realistic values
of $\chi_1, \chi_2, \eta$ in Table~\ref{tab:parameters}~\cite{flavell2014conservative}.
Figure~\ref{fig:CompareCaseNeutral2}
shows the potential $\phi$ 
and 
anion concentration $c_1$.
In this case
with small capacitance parameter $\eta$ = 4.63e-5,
the Robin boundary conditions are almost 
Dirichlet boundary conditions,
so $\phi(\pm 1) \approx -\phi_\pm$
as seen in Fig.~\ref{fig:CompareCaseNeutral2}a.
The anions are attracted to the 
potential maximum at $x = -1$,
leading to a boundary layer in the concentration.
Note that the 
concentration boundary value $c_1(-1) \approx 21.6$
is higher than in case 1 and case 2;
this is attributed to the higher total concentration
$a_1 = 2$ in this case.
These results agree well with the late time results 
for the time-dependent PNP equations in Fig.~2 of~\cite{flavell2014conservative}.

\subsection{Non-electroneutral systems}

This is case 4 in Table~\ref{tab:parameters}.
In contrast to previous cases,
here the total cation concentration $a_2=2$
is greater than the total anion concentration $a_1=1$,
and
the formal potential boundary values are equal,
$\phi_\pm = 1$~\cite{lee2010new}.

\begin{figure*}[htb]
\centering
\includegraphics[width=\textwidth]{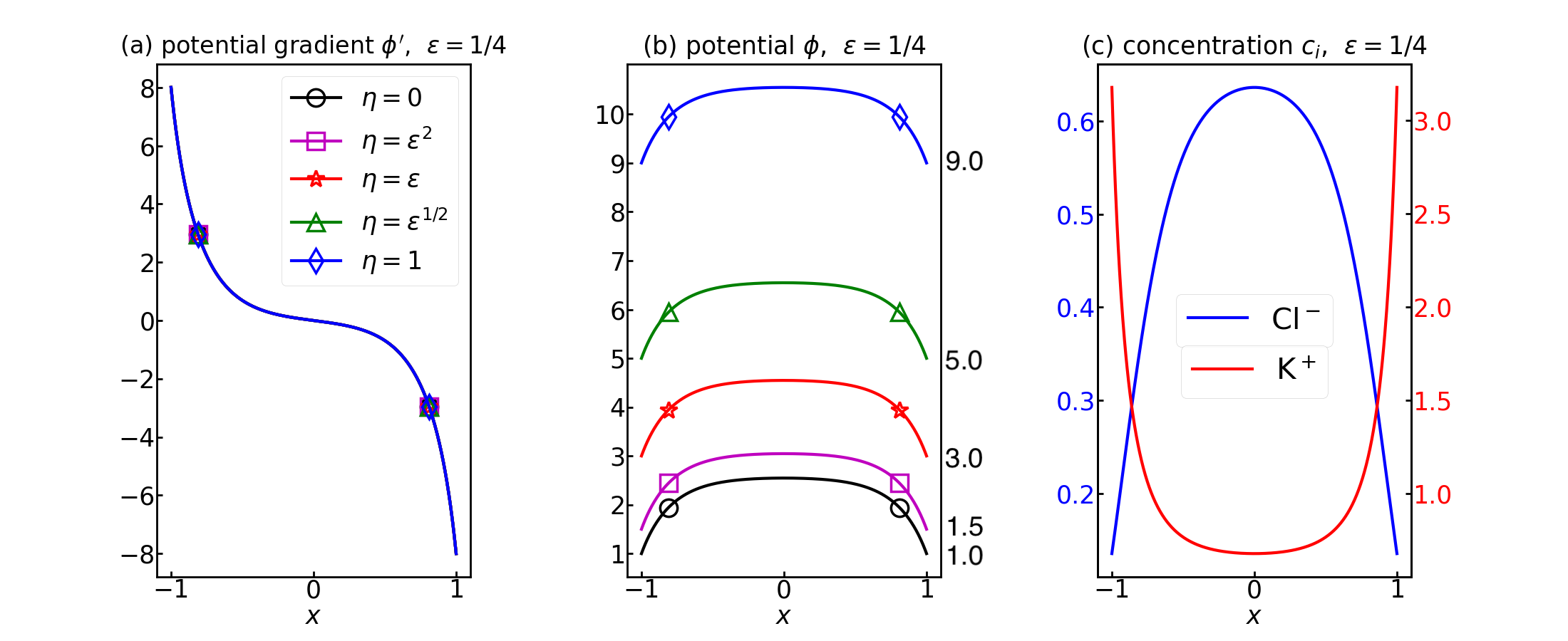}
\includegraphics[width=\textwidth]{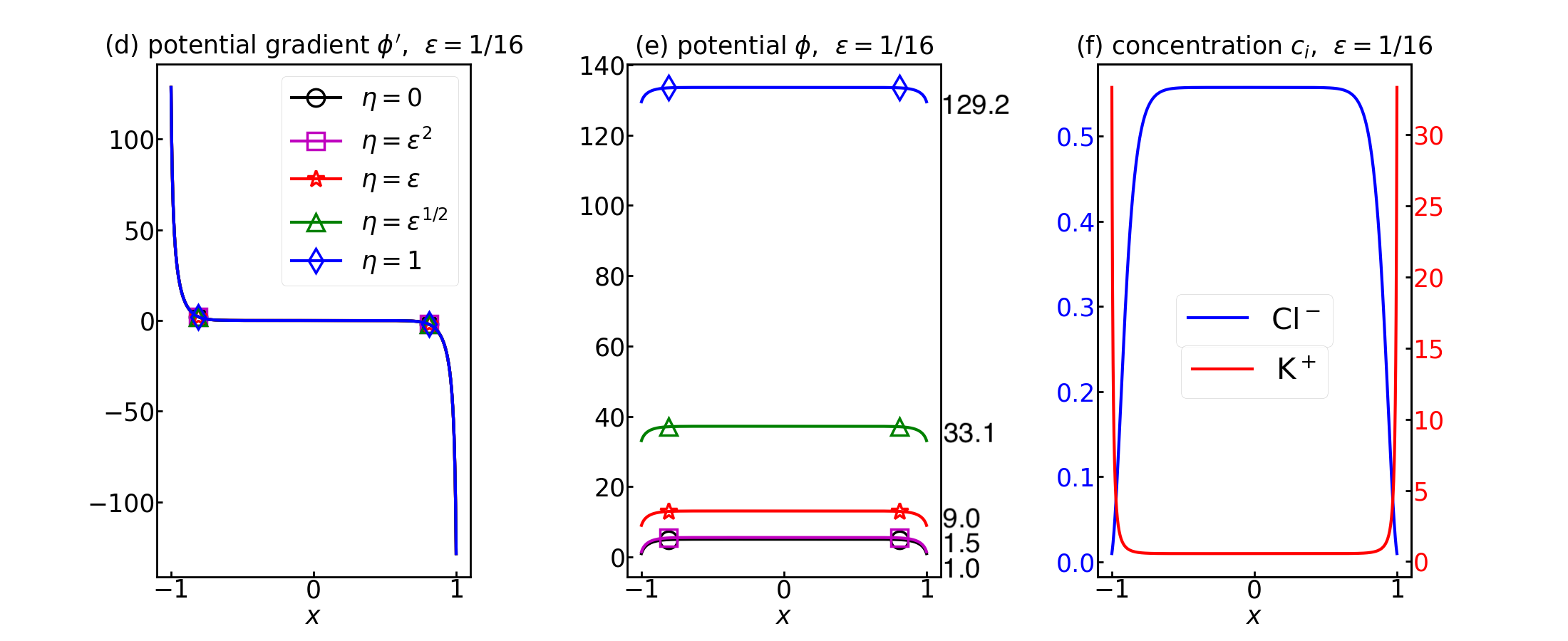}
\caption{
Non-electroneutral systems,
case 4.1, 4.2 in Table~\ref{tab:parameters},
total concentration $a_1=1$ (anion),
$a_2=2$ (cation),
large permittivity $\epsilon = 1/4$ (top), 
small permittivity $\epsilon = 1/16$ (bottom),
equal formal potential boundary values $\phi_\pm = 1$, 
capacitance parameter
$\eta = 0, \epsilon^2, \epsilon, \epsilon^{1/2}, 1$,
(a,d) potential $\phi$,
(b,e) potential gradient $\phi^\prime$,
(c,f) ion concentration 
($c_1$ anion, $c_2$ cation),
integral method with Chebyshev points,
$N=100$ subintervals.}
\label{fig:_non_eletroneutral}
\end{figure*}

Figure~\ref{fig:_non_eletroneutral}
shows the results,
where the top frames use
large permittivity $\epsilon = 1/4$ (case 4.1)
and
the bottom frames use
small permittivity $\epsilon = 1/16$ (case 4.2).
Five values of the capacitance parameter are considered,
$\eta = 0, \epsilon^2, \epsilon, \epsilon^{1/2}, 1$.
Results are shown for the
potential $\phi$,
potential gradient $\phi^\prime$,
and
ion concentrations $c_1, c_2$
(in this case the concentrations have different profiles).
The potential is maximum at $x=0$
and
the potential boundary values $\phi(1)$ 
are printed on the right side of the plots;
these results are in good agreement with~\cite{lee2010new}.
The high potential at $x=0$ 
attracts anions
and 
repels cations,
leading to boundary layers in the ion concentrations.

As the capacitance parameter $\eta$ increases,
the potential also increases,
but the profiles have almost the same shape
and
the potential gradient 
is almost independent of $\eta$
(the curves overlap);
the ion concentrations are also almost
independent of $\eta$
because the 
Nernst-Planck equation~\eqref{eqn:NP_non_dimensional} 
depends on $\phi^\prime$, but not $\phi$.
As the permittivity $\epsilon$ decreases,
the potential gradient
and
ion concentrations develop sharp boundary layers.
In~\cite{lee2010new} 
it was proven that as $\epsilon \to 0$, 
the potential $\phi(x) \to \infty$ for all $x \in [-1,1]$
and
Fig.~\ref{fig:_non_eletroneutral} is consistent with this.
Figure~\ref{fig:_non_eletroneutral} 
also suggests that the 
ion concentration boundary values
converge like 
$c_1(\pm 1) \to 0, c_2(\pm 1) \to \infty$
as $\epsilon \to 0$.

Table~\ref{tab:non_electroneutral}
shows the potential error 
and
convergence rate confirming 2nd order accuracy $p \approx 2$
for this non-electroneutral system.
As before,
for large permittivity $\epsilon = 1/4$
the uniform and Chebyshev points yield similar results.
For small permittivity $\epsilon = 1/16$
the number of iterations is higher for both points sets,
but Chebyshev points
require fewer iterations
and
yield smaller error.

\begin{table}[t!]
\centering
\begin{tabular}{|c|c|c|c|c|c|c|} 
\hline
\multicolumn{7}{|l|}{case 4.1, $\epsilon = 1/4$} \\
\hline
\multicolumn{1}{|c}{$N$} & \multicolumn{3}{|c|}{uniform} & \multicolumn{3}{c|}{Chebyshev} \\
\hline
& iter & error & $p$ & iter & error & $p$ \\
\hline
 50 & 20 & 7.4e-2 &       & 20 & 7.2e-3 & \\ \hline
100 & 20 & 1.8e-2 & 2.039 & 21 & 1.8e-3 & 2.010 \\ \hline
200 & 20 & 4.5e-3 & 2.010 & 21 & 4.5e-3 & 2.005 \\ \hline
400 & 21 & 1.1e-3 & 2.002 & 21 & 1.1e-3 & 2.001 \\ \hline
\hline
\multicolumn{7}{|l|}{case 4.2, $\epsilon = 1/16$} \\ 
\hline
\multicolumn{1}{|c}{$N$} & \multicolumn{3}{|c|}{uniform} & \multicolumn{3}{c|}{Chebyshev} \\
\hline
& iter & error & $p$ & iter & error & $p$ \\
\hline
 50 & nc &        &       & 66 & 1.5e+0  &      \\ \hline
100 & nc &        &       & 64 & 3.9e-1 & 1.940 \\ \hline
200 & 89 & 4.5e-0 &       & 65 & 9.9e-2 & 1.978 \\ \hline
400 & 86 & 8.7e-1 & 2.359 & 68 & 2.3e-2 & 2.088 \\ \hline
\end{tabular}
\vskip 5pt
\caption{
Non-electroneutral system 
(case 4.1, 4.2)~\cite{lee2010new},
parameters in Table~\ref{tab:parameters},
large permittivity $\epsilon = 1/4$, 
small permittivity $\epsilon = 1/16$,
$N$ subintervals, uniform and Chebyshev points,
number of iterations iter, 
potential error $\|{\bf \Phi}_{2N}-{\bf \Phi}_{N}\|_2$, 
convergence rate $p$,
nc = did not converge.}
\label{tab:non_electroneutral}
\end{table}


\section{K$^+$ ion channel}

\begin{figure*}[ht]
\centering
\includegraphics[width=0.8\textwidth]{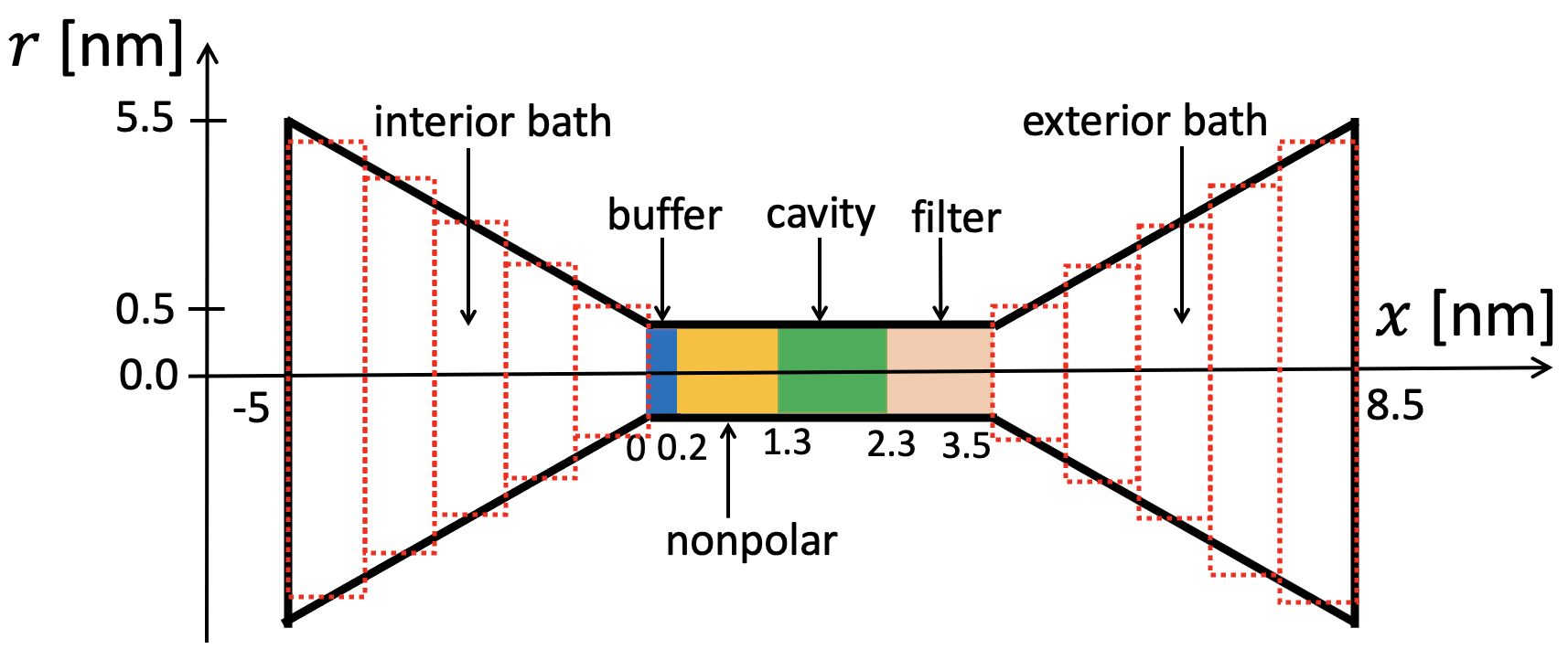}
\caption{
Schematic diagram of 
K$^+$ ion channel model~\cite{gardner2004electrodiffusion}
comprising interior/exterior conical baths with KCl solvent
and
a cylindrical channel with four subregions
(buffer, nonpolar, central cavity, selectivity filter),
channel subregions have constant radius $r = 0.5$ [nm],
bath radius is
linear (black, solid line) in~\cite{gardner2004electrodiffusion}
and
piecewise constant (red, dotted line) in present work.}
\label{fig:KchannelDomain}
\end{figure*}

Finally we consider a potassium (K$^+$) ion channel model
following~\cite{gardner2004electrodiffusion}.
Figure~\ref{fig:KchannelDomain} 
shows a cross-section of the model comprising
two conical baths with KCl solvent
and 
a cylindrical channel with four subregions
(buffer,
nonpolar,
central cavity,
selectivity filter).
In this example the PNP equations are solved
in dimensional form
with the 1D approximation from~\cite{gardner2004electrodiffusion},
\begin{subequations}
\label{PNP_K_channel}
\begin{align}
\label{PNP_K_channel_P}
-&A^{-1}\left(\epsilon A\phi^\prime\right)^\prime = 
e(z_1c_1 + z_2c_2 - \rho_n), \\[2.5pt]
\label{PNP_K_channel_NP}
& A^{-1}\left(A\left(D_ic_i^\prime + 
z_i\mu_ic_i\phi^\prime\right)\right)^\prime = 0, 
\end{align}
\end{subequations}
where $A(x)= \pi r(x)^2$ is the cross-sectional area
in terms of the radius $r(x)$.
Note that the Poisson equation~\eqref{PNP_K_channel_P}
now includes a density of fixed negative charge $\rho_n(x)$.
Table~\ref{tab:K_channel}
gives the parameter values.
The model extends over the interval $[-5, 8.5]$
(the spatial unit is nm).
The four channel subregions have constant radius $r(x) = 0.5$,
and
while~\cite{gardner2004electrodiffusion}
took $r(x)$ to be linear in the baths with $r(-5) = r(8.5) = 5.5$,
here we consider a 
piecewise constant approximation
as suggested by the
red dotted lines in Fig.~\ref{fig:KchannelDomain};
in practice we used 20~steps in each bath.
The other parameters are constant in each subregion,
where the ion mobility is $\mu_i = eD_i/k_BT$
in terms of the diffusion coefficient $D_i$.

\begin{table*}[pt]
\centering
\begin{tabular}{l|c|c|c|c|c|c|c|c|c}
\hline
& interval & $\ell$  & $r(x)$ & $A(x)$   & $Q$ & $\rho_n$ & $\epsilon$ & $\mu_i$ & $D_i$ \\
\hline
interior bath  & [-5.0, 0.0]  & 5 & see text  & $\pi r(x)^2$    & 0     & 0          & 80         & 60    & 1.5 \\
buffer         & [0.0, 0.2]  & 0.2 & 0.5 & $\pi$/4     & -4    & $80/\pi$  & 80         & 16    & 0.4 \\
nonpolar       & [0.2, 1.3] & 1.1 & 0.5  & $\pi$/4     & 0     & 0          & 4          & 16    & 0.4 \\
central cavity & [1.3, 2.3] & 1  & 0.5   & $\pi$/4     & -1/2  & $2/\pi$   & 30         & 16    & 0.4 \\
selectivity filter & [2.3, 3.5] & 1.2 & 0.5 & $\pi$/4     & -3/2  & $5/\pi$   & 30         & 16    & 0.4 \\
exterior bath  & [3.5, 8.5] & 5   & see text  & $\pi r(x)^2$    & 0     & 0          & 80         & 60    & 1.5 \\
\hline
\end{tabular}
\vskip 5pt
\caption{
K$^+$ channel parameters~\cite{gardner2004electrodiffusion},
interval [nm],
length $\ell$ [nm],
radius $r(x)$ [nm],
cross-sectional area $A(x)$ [nm]$^2$,
total fixed charge $Q~[e]$,
fixed charge density $\rho_n$~[molar],
permittivity $\epsilon$,
mobility $\mu_i = eD_i/k_BT$ [1e-5 cm$^2$/(V $\cdot$ s)],
diffusion coefficient $D_i$ [1e-5 cm$^2$/s].}
\label{tab:K_channel}
\end{table*}

The potential and ion concentrations satisfy
Dirichlet boundary conditions,
\begin{subequations}
\label{eqn:K_channel_BC}
\begin{align}
\phi(-5) &= 0, \quad \phi(8.5) = -100 \;{\text{mV}}, \\[2.5pt]
c_i(-5) &= c_i(8.5) = 0.15 \; {\text{molar}}, \quad i = 1,2.
\end{align}
\end{subequations}
There are five interfaces,
$x = 0, 0.2, 1.3, 2.3, 3.5$,
at which the potential 
and ion concentrations satisfy continuity conditions,
\begin{subequations}
\begin{align}
&[\phi] = 0, \quad\left[\epsilon A\phi^\prime\right] = 0, \\[2.5pt]
&[c_i] = 0, \quad 
\left[A\left(D_ic_i^\prime + 
z_i\mu_ic_i\phi^\prime\right)\right] = 0,
\end{align}
\end{subequations}
where the bracket indicates the jump across the interface.

The integral method previously described for systems
with one domain
was extended to this example with
six domains (two baths, four channel subregions)
and
the scheme has the same Gummel iteration form 
as~\eqref{eqn:Gummel}.
In this case the 
computations used uniform grid spacing $h = 0.01$
and
the same stopping criterion as before~\eqref{eqn:stopping_criterion}.
Note that the 
mobility to diffusion coefficient ratio is $\mu_i/D_i = 40$,
which implies that the K$^+$~ion channel problem is 
drift-dominated~\eqref{PNP_K_channel_NP},
so in this case we used a continuation scheme;
the problem was solved for an increasing sequence of ratios,
$\mu_i/D_i = 1, 10, 20, 40$, 
with relaxation parameters $\omega = 0.9, 0.4, 0.26, 0.18$, 
and 
requiring 12, 69, 148, 218 iterations, respectively,
where the converged solution for one ratio $\mu_i/D_i$
is the initial guess for the next ratio;
the total run time for all four ratios was 44 seconds.
The continuation scheme was not needed for the
electroneutral and non-electroneutral systems
considered earlier since the ratio $\mu_i/D_i$
was much smaller in those examples. 

\begin{figure*}[ht!]
\centering
\includegraphics[width=0.8\textwidth]{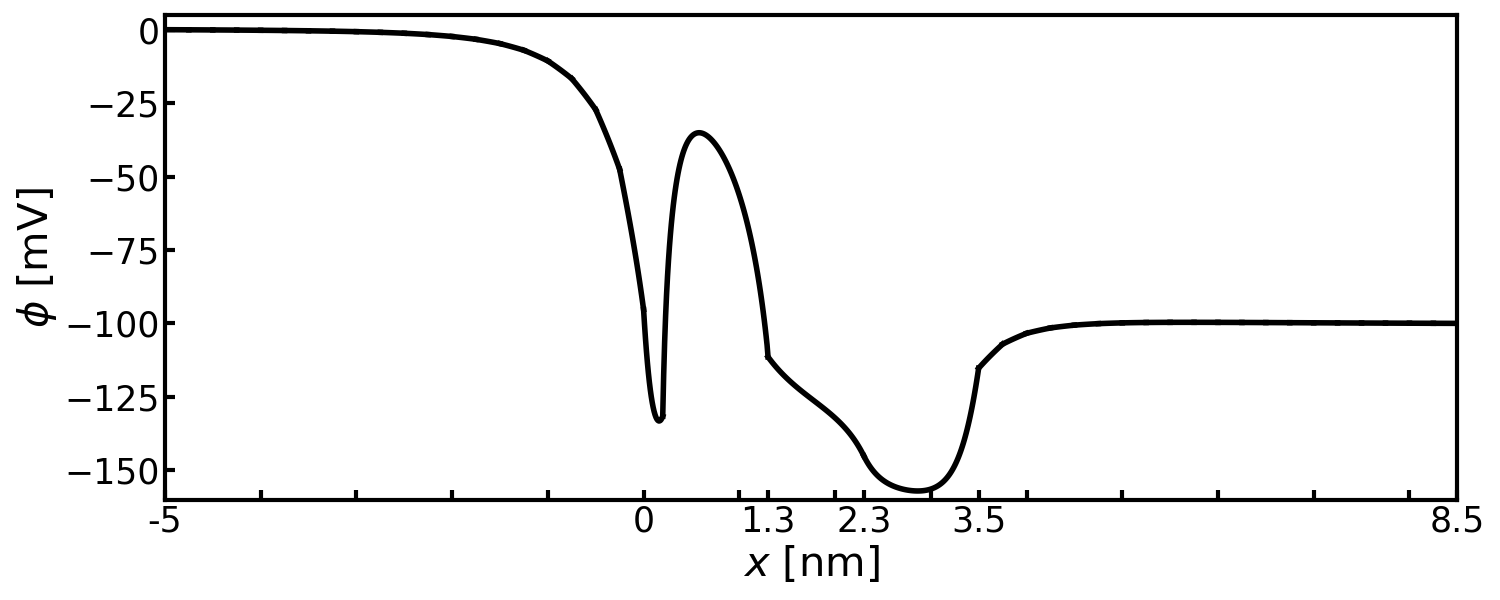}
\includegraphics[width=0.8\textwidth]{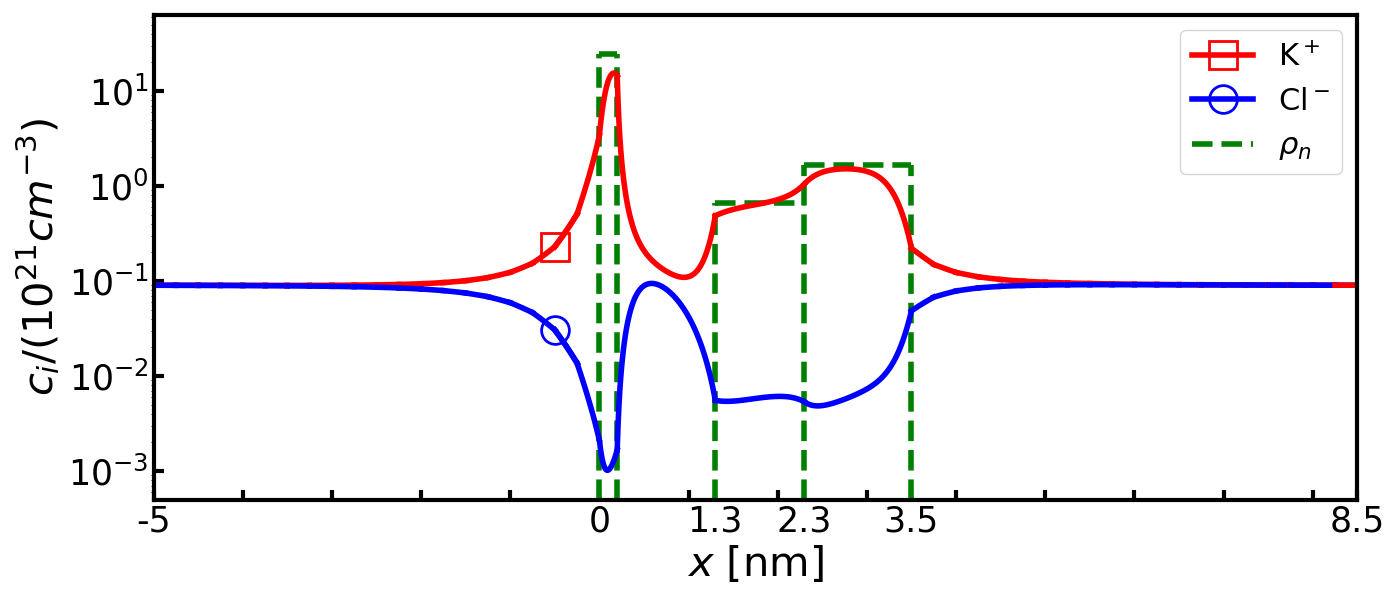}
\caption{
K$^+$~ion channel model 
depicted in Fig.~\ref{fig:KchannelDomain}~\cite{gardner2004electrodiffusion},
parameters in Table~\ref{tab:K_channel},
Dirichlet boundary conditions~\eqref{eqn:K_channel_BC},
fixed negative charge density $\rho_n$
(green dashed lines),
(a) potential $\phi$ [mV],
(b) ion concentration $c_i$ [molar],
K$^+$ (red solid line),
Cl$^-$ (blue solid line),
integral equation method with uniform points, 
grid spacing $h = 0.01$.}
\label{fig:Kchannel}
\end{figure*}

Figure~\ref{fig:Kchannel}
shows the computed potential 
and 
ion concentration profiles
which agree well with the previous 
results~\cite{gardner2004electrodiffusion}.
The negative fixed charge density $\rho_n$,
indicated by green dashed lines,
has substantial impact on the potential
and ion concentrations.
In particular,
the polar subregions with negative fixed charge
have local potential minima below the boundary potential
(-133.7~mV in the buffer, -155.8~mV in the filter),
while the potential rises to a
local maximum in the nonpolar subregion
(-34.3~mV).
As a result,
the K$^+$ cations are attracted to the buffer and filter,
and
are expelled from the nonpolar region,
while the opposite is true for the Cl$^-$ anions;
this demonstrates the cation selectivity property of 
the K$^+$ channel. 

We also computed the current-voltage curve
for the K$^+$~channel model,
where the current is
\begin{subequations}
\begin{align}
\label{current}
I(V_{\rm app}) &= 
-\sum_{i=1}^2 z_ieA(D_ic_i^\prime + z_i\mu_ic_i\phi^\prime), \\
\label{applied_voltage}
V_{\rm app} &= \phi(-5) - \phi(8.5),
\end{align}
\end{subequations}
in terms of the applied voltage $V_{\rm app}$ across the channel.
All other parameters are the same as in Table~\ref{tab:K_channel}.
It follows from the NP equation~\eqref{PNP_K_channel_NP}
that the current~\eqref{current}
is constant throughout the domain
and
in the present calculations this was satisfied to machine precision.

\begin{figure}[ht!]
\centering
\includegraphics[width=0.5\textwidth]{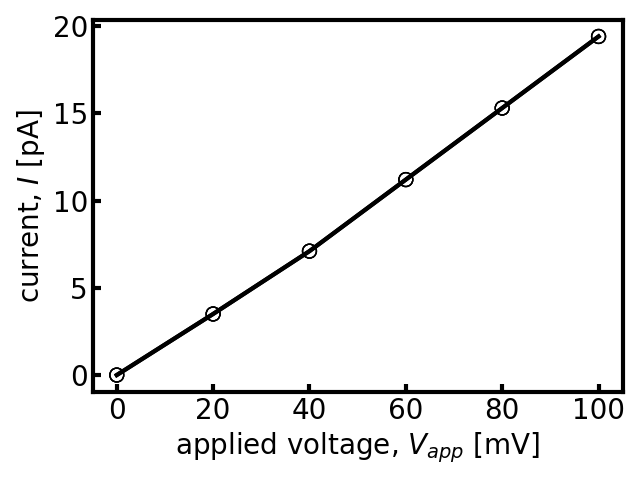}
\caption{
K$^+$~ion channel,
parameters in Table~\ref{tab:K_channel}~\cite{gardner2004electrodiffusion},
current $I$~[pA] versus 
applied voltage $V_{\rm app}$~[mV],
grid spacing $h = 0.01$~nm.}
\label{fig:KchannelIV}
\end{figure}

Figure~\ref{fig:KchannelIV}
shows the computed I-V curve
which closely matches the expected linear relation.
For $V_{\rm app}$ = 100~mV the current is $I$ = 19.4~pA,
which comprises 19.1~pA for K$^+$ cations 
and
0.3~pA for Cl$^-$ anions, 
which again shows the 
cation selectivity of the K$^+$~channel.
For comparison,
\cite{gardner2004electrodiffusion}
obtained 22.2~pA for cations
and
0.3~pA for anions.
In the present case,
we checked that the computed current 
is not sensitive to the number of steps in the bath region,
while Table~\ref{tab:current} shows
that the value approaches the 
result in~\cite{gardner2004electrodiffusion}
as the grid spacing $h$ is reduced.

\begin{table}[hbt]
\centering
\begin{tabular}{c|c|c|c|c|c}
\hline
$h$~nm & 0.02  & 0.01   & 0.005 & 0.0025 & 0.00125 \\
$I$~pA & 14.68 & 19.39 & 20.50 & 20.78 & 20.86 \\
\hline
\end{tabular}
\vskip 5pt
\caption{K$^+$ ion channel model~\cite{gardner2004electrodiffusion},
applied voltage $V_{\rm app}$ = 100~mV,
computed current $I$ [pA] versus grid spacing $h$ [nm].}
\label{tab:current}
\end{table}


\section{Summary}

An integral equation method was presented for the 
1D steady-state PNP equations
modeling ion transport through membrane channels.
Green's 3rd identity was applied to recast
the PNP differential equations as a 
fixed-point problem involving integral operators
for the potential gradient
and
ion concentrations.
The integrals were computed by midpoint and trapezoid rules
at either uniform or Chebyshev grid points,
and
the resulting nonlinear equations were solved by
Gummel iteration with relaxation.
Each iteration requires solving a small linear system for 
the boundary and interface variables.

The integral equation method was applied to
electroneutral and non-electroneutral systems
with Robin boundary conditions on the potential
and
no-flux boundary conditions on the ion concentrations.
The method is 2nd~order accurate
even in cases with sharp boundary layers. 
The method was also applied to a
drift-dominated K$^+$~ion channel model
following~\cite{gardner2004electrodiffusion}.
In this case the fixed charge density
in the channel ensures the 
selectivity of K$^+$~ions
and
exclusion of Cl$^-$ ions.
In these tests the integral method results 
agree well with published
finite-difference and 
finite-element results~\cite{
lee2010new,
flavell2014conservative,
gardner2004electrodiffusion}.

Various techniques are have been used in numerical PNP solvers
to improve stability
such as
positivity-preserving schemes,
Slotboom variables,
and
streamwise upstream differencing.
The integral method calculations presented here 
did not employ these techniques,
although it did use relaxation to accelerate the iteration
and
continuation to compute a good initial guess 
for drift-dominated systems.
It is noteworthy that in these calculations,
the integral method preserved the positivity of the
ion concentrations.

In future work we aim to extend the 
integral equation PNP solver to time-dependent problems
as well as 2D and 3D ion channel models,
where the volume integral can be efficiently computed 
using adaptive quadrature
and
a fast summation method~\cite{
wang2020kernel,
wilson2022tabi,
xu2023dynamics}. 


\bmhead{Acknowledgments}
This work was supported by NSF grants 
DMS-1819094/1819193 and DMS-2110767/2110869.

\section*{Declarations}
\textbf{Conflict of interest} The authors declare that they have no conflict of interest. \\

\noindent\textbf{Data availability} 
The datasets generated 
during the current study can be replicated 
by running the Python scripts in the 
BIE\_PNP\_1D repository, \url{https://github.com/zhen-wwu}.
\bibliographystyle{sn-mathphys}
\bibliography{1DPNP_v2.bib}


\end{document}